\input amstex
\documentstyle{amsppt}
\magnification=1200

\def\ss{\smallskip}
\def\ms{\medskip}
\def\pa{\partial}

\def\O{\Omega}
\def\o{\omega}

\def\crd{\cr\noalign{\vskip4pt}}
\def\eq{\eqalign}

\def\R{${\text{\bf R}}^N$}
\def\no{\noindent}
\def\bs{\bigskip}

\def\g{\gamma}

\def\a{\alpha}
\def\n{\nabla}

\def\da{d_{\partial \O}(x)}
\def\d{\da}

\let\hacek=\v

\nologo
\NoRunningHeads

\topmatter
\title
Comments on a theorem by Ancona
\endtitle
\author 
Andreas Wannebo
\endauthor

\abstract
This note concerns a theorem by A. Ancona, see [1], which gives two
different sufficient conditions on the open set $\O$ in \R\ in order
to make every element in the Sobolev space $W^{m,p}_0(\O)$ a difference
of two nonnegative functions in the space. The proof is consists of two
parts, a main part and an input part. Ancona givs two different inputs
as possible -- a Hardy inequality by Ne\hacek cas and a then new Hardy inequality.

\par
The main part is treated here to give two (new) theorems. 
In this there is no claim of orginality. The scope is improved though.
\par
Newer results by the present author are put into this scheme
and discussed.
\par
A (very) far-reaching conjecture on the main theme introduced by Ancona 
is given.
\endabstract

\endtopmatter

\document

\heading
Introduction
\endheading

This note has as aim to discuss a theorem in a Comte Rendue note by
Alano Ancona, see [1], and most of all its proof, i.e.
what is there and what to do with it?
\par

Two new theorems
are given based on his proof. Since Ancona is brief we are detailed.
There is no claim of orginality on our part here.
The matter is put into the context of later development.
\ss

The note by Ancona had -- already as a manuscript -- 
a big impact on the present author. The note seems to have been 
rather neglected by most authors in the field of Sobolev space theory and 
Hardy inequalities.
\ss

His note was prompted by that Brezis and Browder
needed some positive result in this direction in order to conclude 
their paper, [2], a sequel to a paper on Schr\"odinger operators.
Complications in this second paper depends on the setting in higher order
Sobolev space.
One of the main complications then was the absence of truncation. The 
Ancona note was made to provide a remedy.
\ss

The proof by Ancona involves a general part together 
with Hardy inequalities for domains as input.
He gives two cases for the theorem.
One uses a Hardy inequality for Lipschitz domains by Ne\hacek cas, see
[9], as input and the second one uses a new Hardy inequality that he proved in the note.
\bs

The latter inequality was a major conceptual break through since
no regularity of the boundary is used at all.
At this time it was commonly believed that 
smoothness of the boundary of some kind was needed
for such Hardy inequalities.
\ss

These input results are far from optimal though but the main part 
can be reformulated into a statement, which is optimal in some settings.
\ss

The purpose of this note is mainly to show the importance of the Ancona C.R. note
and to relate to later progress.

\heading
Main part
\endheading

First we fix some notation and defintions.
Let $\O$ denote an open set in \R.
Let $W^{m,p}(\O)$ denote the usual Sobolev space. 
The definition here is that the Sobolev space is the set of all distributions 
with each derivative up to order $m$
equivalent (as distribution) to a real function in $L^p(\O)$.
This is a Banach space with a norm, the Sobolev space norm

$$
||u||_{W^{m,p}(\O)}=\sum_{k=0}^m ||\n^ku||_{L^p(\O)}.
$$
\bs

The Sobolev space $W^{m,p}_0(\O)$ is the closure
of $C^{\infty}_0(\O)$, the infinitely
times differentiable functions with compact support in $\O$,
and the closure is taken in the Sobolev norm,
and thus is a closed subspace of $W^{m,p}(\O)$.
\par

Furthermore let $+$ in the expression 
$$
W^{m,p}_0(\O)_+
$$
denote that the nonnegative cone is taken.

\bs

The following is what Ancona  gives as theorem in [1].
\bs

{\bf 1. Theorem. Ancona.}
\par

{\sl If it holds that $\O$ is open and bounded in \R\ and if it holds that 
$\O$ either is a Lipschitz domain or if $p>N$,
then every $u\in W^{m,p}_0(\O)$
can be written as $u=u_1-u_2$ with $u_i\in W^{m,p}_0(\O)_+$.
\ss

Furthermore there is some constant $c=c(m,p,\O)$ such that

$$
||u_i||_{W^{m,p}(\O)}\le c||u||_{W^{m,p}(\O)}.
$$
}
\ms

The first of these properties can be written 

$$
W^{m,p}_0(\O)=W^{m,p}_0(\O)_+-W^{m,p}_0(\O)_+.
$$
\ms

We will give two theorems by a closer look at the proof by Ancona.
\ms

For more generality we introduce a weighted Sobolev spaces based on a
standard kind of weight functions.
\vfill
\eject

Let $\d=dist(x,\pa\O)$. Then the weighted Sobolev space norm

$$
||u||_{W^{m,p}(\O,\d^s)}
$$

is defined by introducing the weight $\d^s$ into 
all the seminorms that occur as $L^p(\O)$-norms in the definition
the Sobolev space norm.
\ss

However the norm by itself does not determine the elements in a Banach space.
In the unweighted case there is a celibrated
theorem by Meyers and Serrin to the effect that you get the
same Sobolev space $W^{m,p}(\O)$ if you use the definition given above
or use the set of elements defined as the closure of 
$C^{\infty}(\O)$ in the norm.
\ss

The corresponding question for weighted Sobolev space with the norm above
appears not to be settled, i.e. the question if and/or when equality holds.
\ss

This question is avoided here since a subspace is used.
\bs

Define the weighted Sobolev space $W_0^{m,p}(\O,\d^s)$ 
as the closure of  $C_0^{\infty}(\O)$ in the corresponding weighted
Sobolev norm.
\bs

{\bf Definition.} If the supremum of $\d$ is finite when
$x\in \O$ then $\O$ is said to have finite (inner) width.
This is the same as to say that balls inside $\O$ have bounded
radii. 
\bs

Now we turn to the new formulations based on the proof by Ancona.

\bs

{\bf 2. Theorem.}
\par
{\sl
Let $p>1$ and $-\infty < s < \infty $.
Let $\O$ be of finite width,
\par
For 
$u\in W_0^{m,p}(\O,\d^s)$ 
it then holds that
$$
\int_{\O}|u|^p\d^{-mp+s}dx<\infty
$$
implies $u=u_1-u_2$ with
$u_i\in W_0^{m,p}(\O,\d^s)_+$.
}
\bs

This follows together with a norm estimate as in Theorem 1.
\bs

Denote by $W^{m,p}(\O)_{loc}$ the set of functions $\{u\}$ on 
$\O$ with

$$
||u||_{W^{m,p}(\o)}<\infty
$$
for every open $\o$ with $\bar{\o}$ a compact subset of $\O$.
\ss

Then the second theorem reads as follows.
\vfill\eject

{\bf 3. Theorem.}
\par
{\sl Let $p\ge 1$ and $-\infty<s<\infty $. If
$
u\in W^{m,p}(\O)_{loc},
$
then

$$
||u||_{L^p(\O,\d^{-mp+s})}<\infty
$$

and

$$
||\n^mu||_{L^p(\O,\d^s)}<\infty
$$

together implies

$$
u\in W_0^{m,p}(\O,\d^s).
$$
}
\ss

Proofs are given at the end of the paper.

There is some relation here to early papers by Kadlec and Kufner, 
see [4] and [5],
also refered to by Ancona.
\bs

A remark on the history and development of Hardy inequalities
is given in Wannebo [13]. 
It contains a general outline as well as the role of the present author. 
See also Wannebo [10], [11], [12] and [13].
\bs

The main question raised by Ancona is: 
\ms

When holds
$$
W^{m,p}_0(\O)=W^{m,p}_0(\O)_+-W^{m,p}_0(\O)_+?
$$
\ms

Again more generally: 

When holds
$$
W^{m,p}_0(\O,\d^s)=W^{m,p}_0(\O,\d^s)_+-W^{m,p}_0(\O,\d^s)_+?
$$
\bs

One way to get positive answers is to give appropriate Hardy inequalities
to use as inputs to Theorem 2. 
\par

The presently known answers can be found in Wannebo [12].
The results involve conditions given as uniform capacity/uniform 
polynomial capacity conditions.
However a discussion of these concepts is beyond the scope
of this note.
\bs

The following far-reaching conjecture would require more ideas and technique.
\bs

{\bf 4. Conjecture. Wannebo.}
{\sl
Let $\O\subset {\text{\bf R}}^N$ be open.
\ms
(i) For any $m$ odd, any $N$, $p>1$ and $\O$,
it holds

$$
W^{m,p}_0(\O)=W^{m,p}_0(\O)_+-W^{m,p}_0(\O)_+.
$$

(ii) For any $m$ even and positive there exists an $N$, an $\O$ and a $p>1$
such that

$$
W^{m,p}_0(\O)\not= W^{m,p}_0(\O)_+-W^{m,p}_0(\O)_+.
$$
}
\bs

This conjecture is trivial for $m=1$ by truncation.
The case $m=2$ is well-known. 
But this information is too thin.
The conjecture rests instead on theoretical ideas.

The following result is a combination of Theorem 2 and a result
in Wannebo [12].
\bs

{\bf 5. Theorem.}
{\sl
Given a ``certain uniform capacity condition'' on $\pa\O$,
then if $u\in W^{2,p}_0(\O)$, it holds 

$$
u\in W^{2,p}_0(\O)_+-W^{2,p}_0(\O)_+
$$

if and only if

$$
\int_{\O}|u|^p\d^{-2p}dx < \infty.
$$
}

\bs

\heading
Proofs
\endheading

The proofs of Theorem 2 and 3 given below follows from a careful reading 
the proof by Ancona.
\ss

We will need a definition and some notation.
\ss

A Whitney decomposition of $\O$, open, proper subset of \R, 
is a covering of $\O$ with closed cubes. These have pair-wise 
disjoint interior and they shrink as they tend
to the boundery of $\O$ according to the formula ($Q$ any such cube)

$$
diam Q
\le
dist\{Q,\pa \O\}
\le
4diam Q.
$$

We will also use $A$ as a generic constant which is nonnegative
but may vary at each occurence. (Standard.)
\bs

{\sl Proof of Theorem 2.
}
\bs

The goal is to take any function $u\in W^{m,p}_0(\O,\d^s)$
and then to construct a nonnegative majorization with finite
norm. Then the nonpositive part arises as the difference
between these two functions which also has finite norm of course.
\bs

The way to accomplish this is by cutting the original function
into pieces with cut-off functions, to majorize each of these
pieces with a nonnegative funtion and then to sum this new pieces
into a function with finite norm.
\bs

Now we set out to do this procedure.
\bs

Let $Q_0$ be a cube with unit diameter and with centre at the
origin.
\bs

Choose
$\eta\le 0$
with
$
\eta
\in 
C^{\infty}_0
({4\over 3}Q)$
and $\eta|_{Q_0}=1$.
In order to translate to any Whitney cube $Q$,
let $\eta_Q$ be defined as follows
$\eta_Q=\eta(diam(Q)^{-1}(x-x_Q))$ with $x_Q$ centre of $Q$.
Then for $u\in W_0^{m,p}(\O,\d^s)$ denote $u_Q$ for $\eta_Qu$.
\ss

This construction ensures that $u_Q\in W^{m,p}_0({5\over 3}Q)$,
i.e. $u_Q\in W^{m,p}_0(\O,\d^s)$.
\bs

For $p>1$ there is a representation of Sobolev functions
as Bessel potentials. Hence 
$u_Q=G_m*f_Q$ almost everywhere, where $f_Q\in L^p$ and $G_m$
is a certain Bessel kernel. Here it is important that $G_m$
is a nonnegative kernel.
\ss

This two representations of the Sobolev function and they have similar
sized norms, i.e.

$$
||u_Q||_{W^{m,p}}
\sim
||f_Q||_{L^p}.
$$
\ss

Denote $f_{Q,+}=\max\ [f_Q,0]$ and $v_Q=\eta_Q(G_m*f_{Q,+})$.
\ss

This way $v_Q\ge \max \{0, u_Q\}$,
$v_Q\in W_0^{m,p}({5\over 3}Q)_+$,
$v_Q\in W^{m,p}_0(\O,\d^s)_+$.
\ss

In order to avoid messy formulas we denote by $Q'$ the dilation made that
maps $Q$ to a unit cube. When it is clear which cube is refered to,
we write $u'$ etc. for the furter effects of this dilation.
\bs
Now we give a row of inequalities together with the reason that each hold.
The order is the same in both lists. We do the argument first for a unit cube
$Q'$.
\bs
-- Equivalent norms for Sobolev space;
\ss
-- Equivalent norms for Sobolev space and Bessel potentials;
\ss
-- Leibnitz' rule expansion and the triangle inequality;
\ss
-- A Poincar\'e inequality using Bessel potentials;
\ss
-- A triviality for $f_{Q,+}$ and $f_Q$;
\ss
--  The triangle inequality and Poincar\'e inequalities.
\bs

The inequality row is as follows,
$$
\eq{
&||\n^mv_Q'||_{L^p}
\crd
&\sim 
||v_Q'||_{W^{m,p}}
\crd
&\sim
||\n^m\eta_Q'(G_m*f_{Q,+})||_{L^p}
\crd
&\le 
A\sum_{k=0}^m(||\n^k(G_m*f'_{Q,+})||_{L^p}
\crd
&\le
A||f_{Q,+}'||_{L^p}
\crd
&\le
A||f_Q'||_{L^p}
\crd
&\le
A||u_Q'||_{W^{m,p}}
\crd
&\le
A||\n^mu_Q'||_{L^p}
\cr
}
$$
\vfill
\eject

Summing up

$$
||\n^mv'_Q||_{L^p}
\le
A||\n^mu'_Q||_{L^p}.
$$

The expressions here are now dilation homogeneous. Hence we conclude

$$
||\n^mv_Q||_{L^p}
\le
A||\n^mu_Q||_{L^p}.
$$

Next denote

$$
v
=
\sum v_Q,
$$
which is the candidate as the majorizing function. 
\bs

It remains to check the norm.
Since the Sobolev norm is built of some seminorms
it is enough to check them.
\bs

This is done by estimations as a sequence of inequalities.
As before the list of inequalities is preceded by the list of arguments
used.
\bs

Observe that by the very construction and the properties of Whitney cubes
it follows that there is at most a fixed number of overlaps from 
the $\{v_Q\}$.
\bs

List of arguments.
\ss
-- The definition of $v$;
\ss
-- Whitney cube properties;
\ss
-- A Poincar\'e inequality;
\ss
-- Interpolation between Sobolev seminorms;
\ss
-- The bounded inner width of $\O$;
\ss
-- The previous result;
\ss
-- Leibnitz' rule, the triangle inequality finite overlap;
\ss
-- Interpolation cube-wise between seminorms.
\bs

The $p$-power of a seminorm is estimated
\bs

$$
\eq{
&||\n^kv||_{L^p(\O,\d^s)}^p
\crd
&\le
||\sum_Q\n^kv_Q||^p_{L^p(\O,\d^s)}
\crd
&\le
A\sum_Q||\n^kv_Q||^p_{L^p}\cdot l(Q)^s
\crd
&\le
A\sum_Q||\n^mv_Q||^p_{L^p}\cdot l(Q)^{(m-k)p+s}
\crd
&\le
A\sum_Q||\n^mv_Q||^p_{L^p}\cdot l(Q)^s
\crd
&\le
A\sum_Q||\n^mu_Q||^p_{L^p}\cdot l(Q)^s
\crd
&\le
A\sum_{r=0}^m \sum_Q||\n^ru||^p_{L^p}\cdot l(Q)^{-(m-r)p+s}
\crd
&\le
A(||u||_{L^p(\O,\d^{-mp+s})}^p+||\n^mu||_{L^p(\O,\d^s)}^p).
\cr}
$$
\bs

Now the first expression is finite according to the assumption
on $u$ and the second is finite since $u\in W^{m,p}(\O,\d^s)$.
\bs

End of proof.
\bs

{\sl
Proof of Theorem 3.
}

The proof is quite the same. Just take $u\in W^{m,p}(\O)_{loc}$ instead
and then proceed as before. The summation of the $v_Q$ then gives a result
in the closure of $C^{infty}_0(\O)$ by the convergence and with the right
norm because of the estimates.

Furthermore the proceedure can be simplified since nonnegativity and
the Bessel kernel is not needed. The argument can be made solely by usual
seminorms in Sobolev space.
\bs

End of proof.
\bs
\bs

For the record we give a list of misprints in Wannebo [10]. The list was given in Wannebo 1991 thesis: ``Some topics in Sobolev space theory''.
(All items have been publicized.)
\bs
p.90 line -10 ball has as centre the centre of $Q'$;
\par
p.90 line -17 reads $x$ in the interior;
\par
p.90 line -7 reads $k+a$;
\par
p.91 line -7 reads $\int_Q|\n^k(u-P)|^pdx$ with $P$ polynomial of 
degree $\le m-1$;
\par 
p.91 line -1 reads $2^{(m-1)p}$;
\par
p.92 line -12 reads $\sum_{n=a}^{s+a}$;
\par
p.92 line -9 reads with extra factor $e^{a\a}$ (-- no consequence);
\par
p.92 line 10 reads $(\varrho(x)-n+a+1)$ in integrand;
\par
p.93 line 16 reads $|\g|$;
\par
p.93 line 93 line -6 reads $\sum_{r=0}^{k-1}$;
\par
p.94 line -5 reads $A\cdot\g_{m,m-1,p}(K,2Q)$.
\bs

\enddocument
\bye